\documentclass[11pt, a4paper]{article}

\usepackage{amsmath,amssymb,amstext}
\usepackage{graphicx,color}
\usepackage{graphics}

\begin{document}

\title{Probabilistic Analysis of LCF Crack Initiation Life of a Turbine Blade under Thermomechanical Loading}

\author{S. Schmitz$^{\footnotesize\textsf{1, 2}}$, G. Rollmann$^{\footnotesize\textsf{1}}$, H. Gottschalk$^{\footnotesize\textsf{3}}$, R. Krause$^{\footnotesize\textsf{2}}$ \\
        $^{\footnotesize\textsf{1}}$ Siemens AG, Mellinghoferstra\ss e 55,\\ 45473 M\"ulheim an der Ruhr,
        Germany\\
        $^{\footnotesize\textsf{2}}$ Universita della Svizerra Italiana, \\Institute of Computational Science,\\
        Via Giuseppe Buffi 13, 6904 Lugano, Switzerland\\
        $^{\footnotesize\textsf{3}}$ Bergische Universit\"at Wuppertal,\\ Department of Mathematics and Computer Science,\\
        Gau\ss stra\ss e 20, 42097 Wuppertal, Germany}

\maketitle


\begin{abstract}
An accurate assessment for fatigue damage as a function of
activation and deactivation cycles is vital for the design of many
engineering parts. In this paper we extend the probabilistic and
local approach to this problem proposed in \cite{Schmitz_Seibel},
\cite{Gottschalk_Schmitz} and \cite{Schmitz_Rollmann} to the case
of non-constant temperature fields and thermomechanical loading.
The method has been implemented as a finite element postprocessor
and applied to an example case of a gas-turbine blade which is
made of a conventionally cast nickel base superalloy.



\end{abstract}


\section{Introduction}
\label{section_introduction}

The necessity for a flexible service of a lot of engineering parts
such as gas turbines leads to the importance of fatigue analysis,
where probabilistic models can be very valuable.
In this work, we present a probabilistic model for low-cycle
fatigue (LCF) which can be derived from the Poisson point process
or from a spatial hazard approach, confer \cite{Schmitz_Seibel} and
\cite{Gottschalk_Schmitz}, respectively.
The model can be applied to polycrystalline metal which is
sufficiently fine-grained so that isotropic material behavior can
be assumed and continuum mechanics can be employed. Here, failure
of a component is defined to be given by the initiation of the first LCF crack.
The probabilistic model yields the probability of failure
(PoF) as a function of the number of load cycles.
We have extended the model proposed in \cite{Schmitz_Seibel},
\cite{Gottschalk_Schmitz} and \cite{Schmitz_Rollmann} by a
temperature model of the LCF parameters as well as the percentile
bootstrap method in order to be able to consider uncertainties due
to LCF test data. In contrast to the deterministic safe-life
approach \cite{Harders_Roesler}, this extended probabilistic LCF
model takes inhomogeneous temperature and strain fields, size
effects and uncertainties due to specific calibration data into
account.

We apply the probabilistic model to a gas-turbine blade which is
subjected to thermomechanical loading during the operating state.
In this case, the corresponding LCF failure
mechanism is surface driven. Having computed
the total PoF we also consider and visualize the hazard density on
the blade's surface.

\section{A Probabilistic Model for LCF}
\label{section_probabilistic_LCF}

In the following, we first briefly revisit the spatial
hazard approach for surface driven LCF presented in
\cite{Gottschalk_Schmitz} and \cite{Schmitz_Rollmann}
and consider LCF as a
failure-time process. If $N$ denotes the random variable which
represents the cycle of first crack initiation and $P$ the
underlying probability measure, the hazard rate is defined by
\begin{equation}\label{probabilistic.1.1}
    h(n)=\lim_{\Delta n\rightarrow0}\frac{P(n<N\leq n+\Delta
n|N>n)}{\Delta n}=\frac{f_N(n)}{1-F_N(n)},
\end{equation}
confer \cite{Escobar_Meeker}. Here, we model $N$ as a continuous
random variable in agreement with the literature, confer \cite{Harders_Roesler}
and \cite{Vormwald}. $F_N(n)=P(N\leq n)$ is
the cumulative distribution function and $f_N(n)=dF_N(n)/dn$ the corresponding
density function. The hazard rate $h$ is also called instantaneous
failure rate as for a small step $\Delta n$ the expression
$h(n)\cdot\Delta n$ is an approximation for the propensity of
failure in the next time step $\Delta n$, given no failure to time
$n$. 

Considering strain controlled LCF failure mechanism on a component
which is made of polycrystalline metal and represented by a domain
$\Omega$, we assume according to \cite{Gottschalk_Schmitz} and \cite{Schmitz_Rollmann}
that the surface zone which is affected from
the crack initiation process of a single LCF crack is small with
respect to the surface of the component. Thus, we suppose that
in any subregion $A$ of the component's surface
$\partial\Omega$, the corresponding hazard rate $h_A$ is a local
functional of the displacement field $\textbf{u}$ and the
temperature field $T$ in that particular region with
\begin{equation}\label{probabilistic.2.1}
h_{A}(n)=\int_{A} \rho(n;\nabla \textbf{u},T)\, dA.
\end{equation}
Here, $\nabla \textbf{u}$ is the Jacobian matrix of $\textbf{u}$.
We call the integrand $\rho$ hazard density function which is the
core of this spatial hazard approach.
For inhomogeneous strain fields
$\varepsilon_a=\varepsilon_a(\nabla \textbf{u},T)$ we obtain
$h(n)=\int_{\partial\Omega}\rho(n;\varepsilon_a,T)\,dA$. Here,
$\varepsilon_a$ is an equivalent strain amplitude which can be
derived from thermoelastic finite element analysis (FEA) with
subsequent application of stress-strain relationships. For more
details confer \cite{Schmitz_Rollmann}, \cite{Harders_Roesler} and
\cite{Neuber}.


Taking $F_N(n)=1-\exp\left(-\int_0^nh(s)\,ds\right)$ into account
the probability of LCF crack initiation on the surface
$\partial\Omega$ until cycle $n$ is given by
\begin{equation}\label{probabilistic.3.1}
F_N(n)=1-\exp\left(-\int_0^n\int_{\partial\Omega}\rho(s;\varepsilon_a,T)\,dAds\right).
\end{equation}
In Section \ref{section_turbine_blade} the hazard density function $\rho$ will
be employed to identify the critical and possibly overengineered regions of the component.
Considering the statistical evalutaions in \cite{Schmitz_Seibel} and
following \cite{Gottschalk_Schmitz} and \cite{Fedelich} we
assume that the number $N$ of cycles to crack initiation are
Weibull distributed.
Therefore, we choose the Weibull hazard ansatz
\begin{equation}\label{probabilistic.4.1}
\rho(n;\varepsilon_a,T)=\frac{m}{N_{\textrm{det}}(\varepsilon_a,T)}\left(\frac{n}{N_{\textrm{det}}(\varepsilon_a,T)}\right)^{m-1},
\end{equation}
where $m$ is the Weibull shape and $N_{\textrm{det}}(\varepsilon_a,T)$
determines the corresponding Weibull scale parameter.
If the Coffin-Manson-Basquin (CMB) equation
-- confer \cite{Harders_Roesler} and \cite{Vormwald}  --
is taken for the strain-life relationship the scale field
$N_{\textrm{det}}(\mathbf{x})=N_{\textrm{det}}(\varepsilon_a(\mathbf{x}),T(\mathbf{x}))$
is given by the solution of
\begin{equation}\label{probabilistic.5.1}
\varepsilon_a(\mathbf{x})=\frac{\sigma_f'(T(\mathbf{x}))}{E(T(\mathbf{x}))}(2N_{\textrm{det}}(\mathbf{x}))^{b(T(\mathbf{x}))}+
\varepsilon_f'(T(\mathbf{x}))(2N_{\textrm{det}}(\mathbf{x}))^{c(T(\mathbf{x}))}
\end{equation}
on every point $\mathbf{x}\in\partial\Omega$, where the surface
$\partial\Omega$ is subjected to an equivalent strain field
$\varepsilon_a(\mathbf{x})$ and a temperature field
$T(\mathbf{x})$. Having chosen an appropriate temperature
model\footnote{In this work we use a proprietary temperature model
by Siemens AG.} for the CMB parameters
$\sigma_f',b,\varepsilon_f',c$ and for Young's modulus $E$,
the probabilistic model for LCF is given by the Weibull
distribution
\begin{equation}\label{probabilistic.6.1}
F_N(n)=1-\exp\left[-\left(\frac{n}{\eta}\right)^m\right]\quad\textrm{for
scale}\quad
\eta=\left(\int_{\partial\Omega}\frac{1}{N_{\textrm{det}}^{\,\,
m}}dA\right)^{-1/m}
\end{equation}
and for some shape parameter $m\geq1$, which yields the
probability for LCF crack initiation to cycle $n$.


The model can be calibrated by means of usual maximum likelihood
methods, confer \cite{Escobar_Meeker} and \cite{Georgii}, for
example. The parameter $m$ determines the scatter of the
distribution where small values for $m\geq1$ correspond to a large
scatter. 

The CMB parameters of the model are not the same as obtained from
fitting standard specimen data. Due to the size effect, see
\cite{Vormwald}, the original CMB approach leads to different
values of the parameters for specimens under the same temperature
and strain conditions but with different gauge areas. The new and
more physical interpretation of these parameters -- in context of
the probabilistic model according to \cite{Schmitz_Seibel}  and
\cite{Gottschalk_Schmitz} -- takes the size effect and inhomogeneous
strain and temperature fields into account, so that the CMB
parameters can be calibrated with LCF-test results of specimens
with arbitrary geometry\footnote{The specimen must consist of
sufficiently many grains so that continuum
mechanics can be applied. Moreover, information on when the first
LCF crack initiation occurred has to be provided which can be
practically difficult, however.} and under arbitrary strain and
temperature fields. Thus, these newly interpreted parameters can
be assigned to every such geometry. 

\section{LCF Crack Initiation Life of a Turbine Blade}
\label{section_turbine_blade}

In this section we consider a turbine blade which is made of a
polycrystalline cast nickel-based superalloy such as RENE 80. The
following probabilistic analysis of its LCF crack initiation life
is based on an FEA model for the operating state and on the
Weibull distribution (\ref{probabilistic.6.1}).

\begin{figure}[htbp]
  \centering
    \includegraphics[width=4.0cm]{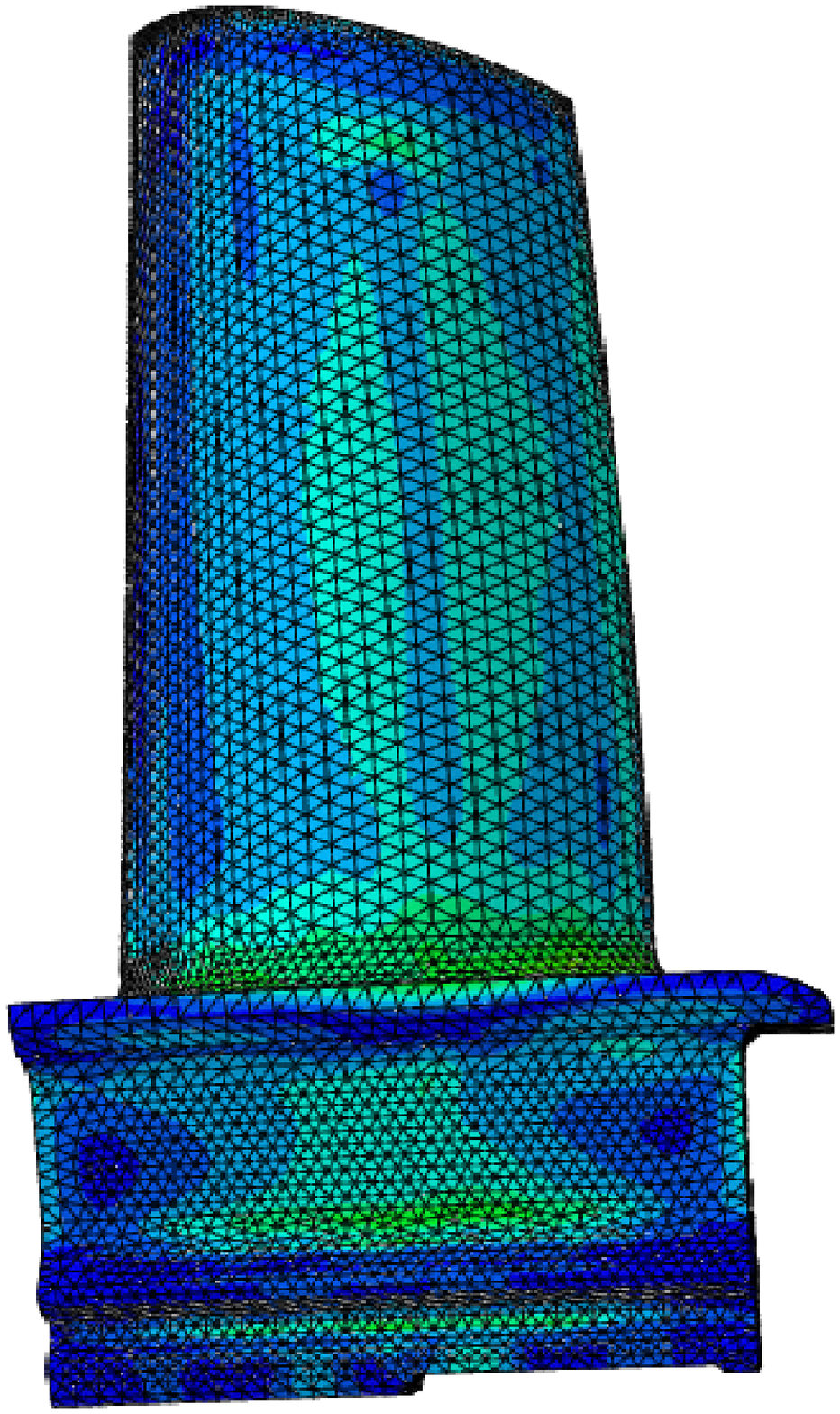}   
    \hspace{0.5cm}
    \includegraphics[width=4.0cm]{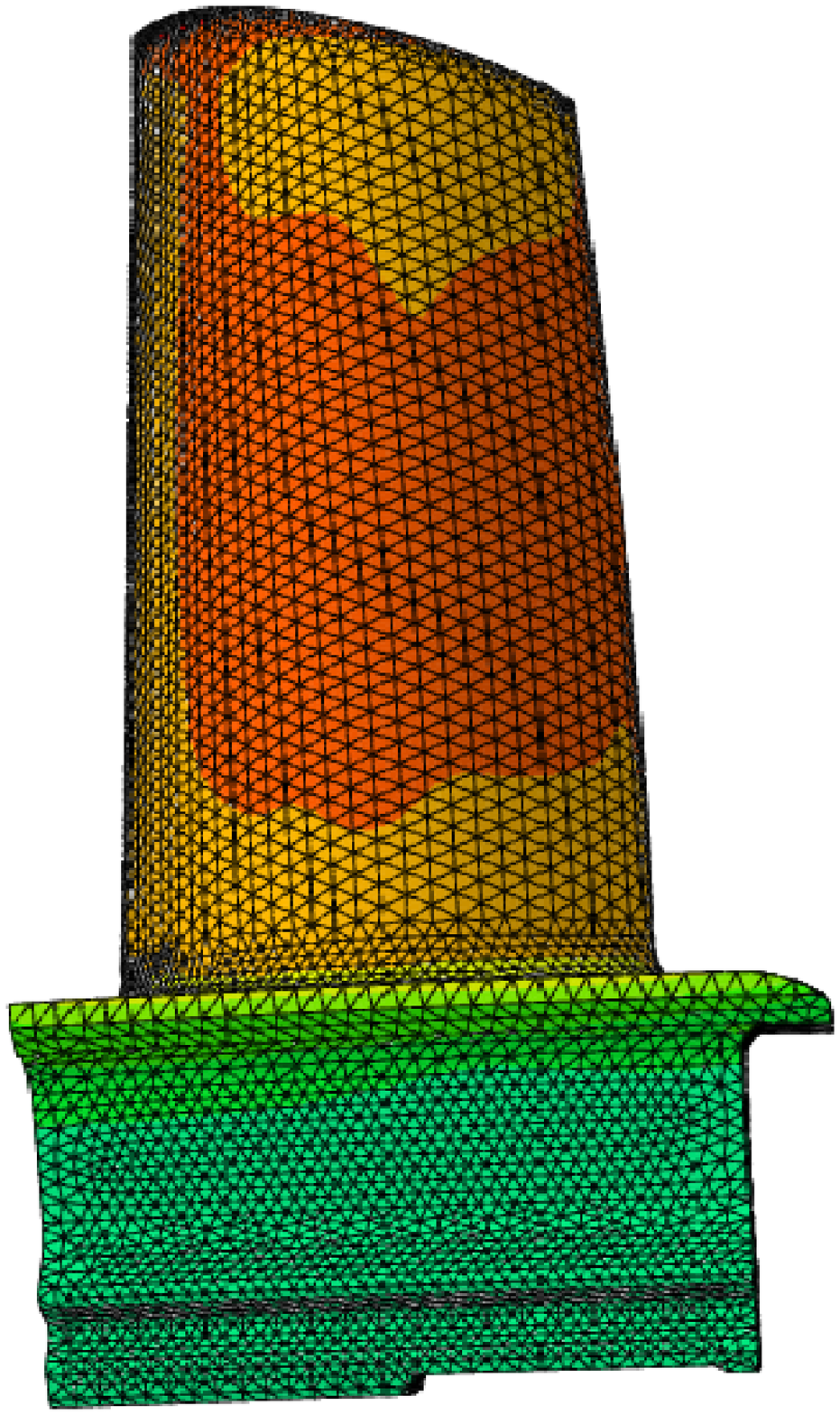}   
  \caption{FEA results of Abaqus 6.9-2 for the von Mises stress field (left)
and temperature field (right) of the turbine blade.}
  \label{Figure_FEA}
\end{figure}

The temperature and strain field of the turbine blade in the
operating state is computed by means of a thermoelastic FEA-model
within Abaqus 6.9-2 and of the Neuber shakedown
method\footnote{Confer \cite{Harders_Roesler} and \cite{Neuber}.}
which considers plasticity. The model includes approximately
190,000 tetrahedral, affine Lagrange elements of Serendipity class
with 10 nodes. Figure \ref{Figure_FEA} shows the von Mises stress
and temperature field in the operating state. In the shutdown
state the von Mises stress is everywhere zero and the temperature
field is set equal to that one of the operating state which is an
conservative approximation and avoids the treatment of
thermo-mechanical fatigue (TMF). The transition from the shutdown
state to the operating state and then back to the shutdown state
is considered as one load cycle. It is further assumed that the
shutdown and operating state stay the same during the cycles.


The probabilistic LCF model has been calibrated with strain
controlled LCF test results\footnote{The results were provided by
Siemens AG.} for standard specimens which were subjected to
different temperatures and strain amplitudes. The maximum
likelihood method has been used for the calibration which
constitutes a statistical estimation of the model parameters,
confer \cite{Escobar_Meeker}. This estimation in conjunction
with computing the surface integral (\ref{probabilistic.6.1})
results in a value for the Weibull scale parameter $\eta$ which
we call the maximum likelihood value for $\eta$.

As the estimation depends on the LCF
test results there are uncertainties for the values of the model
parameters. This affects the total PoF with respect to LCF crack
initiation. We employ the fully parametric bootstrap sampling
procedure in conjunction with the percentile method -- confer
\cite{Escobar_Meeker} -- to consider theses additional
uncertainties. We used 2,000  bootstrap samples which were
obtained from the maximum likelihood estimation.
These samples are
different parameter realizations of the probabilistic model for
LCF and describe the distributions of the model parameters.
Using these parameter realizations in conjunction with the
FEA postprocessing results in 2,000 values for the Weibull shape
and scale parameter. Then, the law of total probability
yields the total PoF with respect
to LCF crack initiation. Note that computationally most expensive parts of the FEA
postprocessing are identifying the blade's surface and computing
the surface integral in (\ref{probabilistic.6.1}) for the
bootstrap samples.


\begin{figure}[t]
    \centering
\scalebox{0.32}{\includegraphics{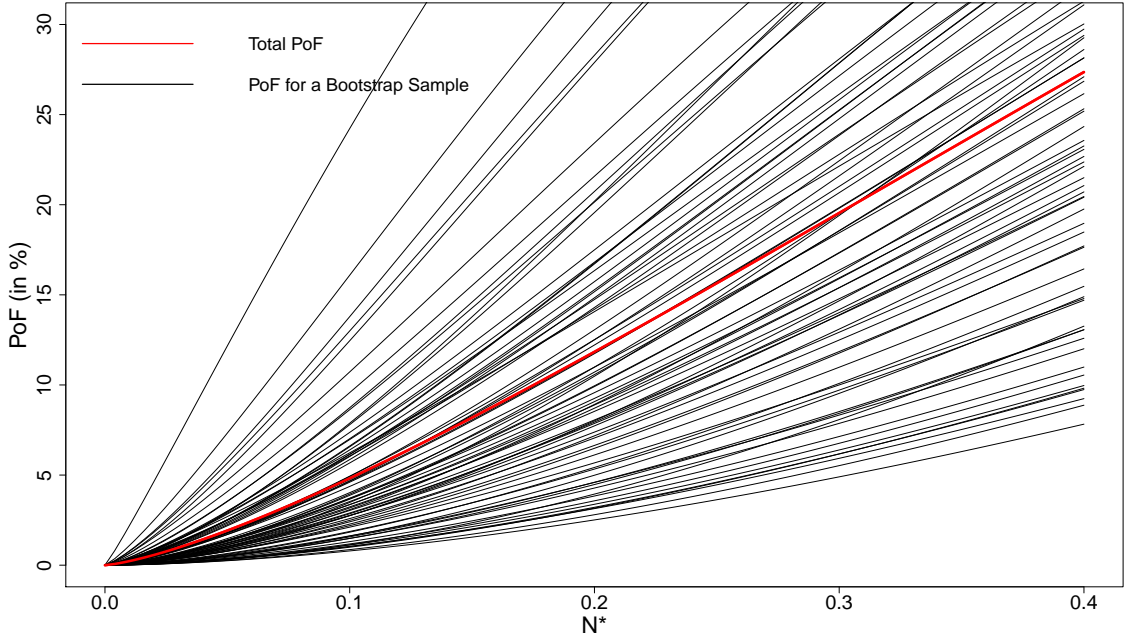}}
\caption{Total PoF (red) due to LCF crack initiation and single
PoF curves (cumulative Weibull distributions, black) corresponding
to 72 of the 2,000 bootstrap samples.
$N^*$ are multiples of the maximum likelihood value for the Weibull scale $\eta$.} 
    \label{Figure_PoF}
\end{figure}

\begin{figure}[t]
    \centering
\scalebox{0.55}{\includegraphics{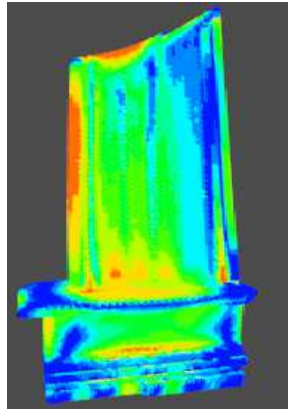}}  
\caption{Hazard density of the analyzed turbine blade: Red areas
show critical regions and blue ones may indicate overengineered
regions of the blade.}
    \label{Figure_hazard}
\end{figure}

Figure \ref{Figure_PoF} shows cumulative Weibull distributions
corresponding to 72 of the 2,000 bootstrap samples. Each black
curve is one parameter realization of the probabilistic model and
yields different values for the PoF depending on multiples $N^*$
of the maximum likelihood value of $\eta$. The law of total probability
results in the total PoF (red curve).
For $N^*=0.0796$ the total PoF is $3.032\%$, for 
example. From a design perspective one decides which PoF is
acceptable and then chooses the corresponding number of allowable
shutdown and operating cycles. Figure \ref{Figure_hazard} shows
the hazard density on the turbine's surface, where the red regions
are areas with higher risk for crack initiation. Blue regions in some cases
may indicate overengineered areas of the turbine blade. To obtain a
final judgment on overengineered areas a more complete design
perspective is needed which includes performance and efficiency
criteria, for example.



\section{Conclusion}
\label{section_conclusion}

In this work, we computed the hazard density and the total PoF of
a turbine blade under cyclic loading due to LCF crack initiation
on the surface. The computations are based on the probabilistic
model for LCF as proposed in \cite{Schmitz_Seibel},
\cite{Gottschalk_Schmitz} and \cite{Schmitz_Rollmann}. In order to
consider thermomechanical loading and uncertainties due to LCF
test data we extended the model by a temperature model and
included the percentile bootstrap method, respectively.

In future, we plan to consider information on
local strain gradients. LCF cracks initiating at the surface grow
into the component, where a different local strain field may
result to a different speed of crack growth, confer
\cite{Harders_Roesler}. Moreover, we plan to extend the probabilistic model to
consider HCF, TMF and non-stationary FEA. Finally, note that the model can be
used to optimize the total PoF with respect to the shape $\Omega$, i.e. to
find a design $\Omega$ under certain constraints such that surface
integrals of the form of (\ref{probabilistic.6.1}) are minimized.
This is also called
optimal reliability, confer \cite{Gottschalk_Schmitz}. 




\subsubsection*{Acknowledgment}
This work has been supported by the German federal ministry of
economic affairs BMWi via an AG Turbo grant. We wish to thank the
gasturbine technology department of the Siemens AG and the
institute of energy and climate research (IEK-2) of the FZ J\"ulich
for stimulating discussions and many helpful suggestions.





\end{document}